\documentclass[12pt]{article}
\usepackage{amsmath}
\usepackage{epsf}
\usepackage{amssymb}
\usepackage{graphicx}
\usepackage{ifpdf}

\newtheorem{theorem}{Theorem}
\newtheorem{definition}[theorem]{Definition}
\newtheorem{lemma}[theorem]{Lemma}
\newtheorem{observation}[theorem]{Observation}
\newtheorem{claim}[theorem]{Claim}

\begin{document}
\pagestyle{empty}
\renewcommand{\thefootnote}{\fnsymbol{footnote}}

\begin{center}
{\bf \Large The crossing number of the generalized Petersen graph $P(3k,k)$ in the projective plane}\footnote{This
work was supported by Hunan Provincial Natural Science Foundation(No. 2018JJ2454) and Hunan Education Department Foundation(No. 18A382).}
\vskip 5mm

{{\bf Jing Wang,   Zuozheng Zhang\footnote{Correspondence author}}\\[2mm]
  Department of Mathematics and Computer Science,
Changsha University, Changsha 410003, China}\\[6mm]
\end{center}
\date{}

\noindent{\bf Abstract} The crossing number of a graph $G$ in a surface $\Sigma$, denoted by $cr_{\Sigma}(G)$, is the minimum number of pairwise intersections of edges in a drawing of $G$ in $\Sigma$. Let $k$ be an integer satisfying $k\geq 3$, the generalized Petersen graph $P(3k,k)$ is the graph with vertex set $V(P(3k,k))=\{u_i, v_i| i=1,2,\cdots,3k\}$ and edge set $E(P(3k,k))=\{u_iu_{i+1}, u_iv_i, v_iv_{k+i}| i=1,2,\cdots,3k\},$
the subscripts are read modulo $3k.$ This paper investigates the crossing number of $P(3k,k)$ in the projective plane. We determine the exact value of $cr_{N_1}(P(3k,k))$ is $k-2$ when $3\le k\le 7,$ moreover, for $k\ge 8,$ we get that $k-2\le cr_{N_1}(P(3k,k))\le k-1.$

{\bf MR(2000) Subject Classification} 05C10, 05C62

\section{Introduction}
\label{secintro}
A {\it surface} $\Sigma$ means a compact, connected 2-manifold. It is known that there are two kinds of closed surfaces, orientable and nonorientable \cite{JG1987}. Every closed connected orientable surface is homeomorphic to one of the standard surfaces $S_k$ $(k\ge 0)$, while each closed connected nonorientable surface is homeomorphic to one of $N_k$ $(k\ge 1)$. In particular, the {\it projective plane}, $N_1$, is a 2-manifold obtained by identifying every point of the 2-sphere with its antipodal point.

Let $G=(V,E)$ be a simple graph with vertex set $V$ and edge set $E$. Let $D$ be a good drawing of the graph $G$ in a surface $\Sigma,$ we denote the number of pairwise intersections of edges in $D$ by $v_D(G:\Sigma)$, or by $v(D)$ without ambiguous. The {\it crossing number of $G$ in a surface $\Sigma$}, denoted by $cr_{\Sigma}(G)$, is the minimum number of pairwise intersections of edges in a drawing of $G$ in $\Sigma$, i.e.,
$$cr_{\Sigma}(G)=\min_{D}v_D\big(G:\Sigma\big).$$
In particular, the {\it crossing number of $G$ in the plane $S_0$} is denoted by $cr(G)$ for simplicity. It is well known that the crossing number of a graph in a surface $\Sigma$ is attained only in {\it good drawings} of the graph, which are the drawings where no edge crosses itself, no adjacent edges cross each other, no two edges intersect more than once, and no three edges have a common point.

In a drawing $D$ of $G=(V,E)$ in a surface $\Sigma,$ if an edge is not crossed by any other edge, we say that it is {\it clean} in $D,$ otherwise, we say it is {\it crossed}. Moreover, let $F \subseteq E$, we say $F$ is clean in $D$ if all of the edges in $F$ are clean, otherwise, we say  $F$ is crossed.

Let $A$ and $B$ be two (not necessary disjoint) subsets of the edge set $E,$ the number of crossings involving an edge in $A$ and another edge in $B$ is denoted by $v_D(A,B).$ In particular, $v_D(A,A)$ is denoted by $v_D(A).$ By counting the number of crossings in $D,$ we have

\begin{lemma}\label{leAB}
Let $A,B,C$ be mutually disjoint subsets of $E,$ then
\begin{eqnarray*}
v_D(A,B\cup C)= v_D(A,B)+v_D(A, C),\\
v_D(A\cup B)= v_D(A)+v_D(A,B)+v_D(B).
\end{eqnarray*}
\end{lemma}

Computing the crossing number of a given graph is, in general, an elusive problem. Garey and Johnson have proved that the problem of determining the crossing number of an arbitrary graph in the plane is NP-complete \cite{GM1983}. Because of its difficulty, there are limited results concern on this problem, see \cite{Kleitman1970,MK2001,Bokal12007,LX2005,MDJ2005} and the references therein. The generalized Petersen graph $P(n,k)$ is a counterexample to many conjectures and thus plays an important role in graph theory. Exoo, Harary and Kabell began to study the crossing number of $P(n,k)$ in $S_0$ and they worked out the case when $k=2$ \cite{GE1981}. For the case $k=3$, Fiorini proved that $cr(P(9,3))=2$ \cite{FS1986}, later on, Richter and Salazar \cite{RB2002} determined the crossing number of $P(3t+h,3)$ is $t+h$ if $h\in\{0,2\}$ and $t+3$ if $h=1$, for each $t\ge 3$, with the single exception of $P(9,3)$. We tried to obtain $cr(P(n,k))$ when $n$ can be expressed as a function of $k$, and proved that $cr(P(3k,k))=k$ for $k\ge 4$ \cite{Wang2011}.

As for the crossing number of graphs in a surface other than the plane, it is not surprising that the results are even more restricted: only the crossing number of Cartesian product graph $C_3\square C_n$ in $N_1$ \cite{Riskin1993}, the crossing number of the complete graph $K_9$ in $S_2$ \cite{Riskin1995}, the crossing number of the complete bipartite graph $K_{3,n}$ in a surface with arbitrary genus \cite{K3n}, the crossing number of $K_{4,n}$ either in $N_1$ or in $S_1$ \cite{PTH2005,PTH2009}, and the crossing number of  the circulant graph $C(3k; \{1,k \})$ in $N_1$ \cite{C3k} have been determined. These facts motivate us to investigate the crossing number of $P(3k,k)$ ($k\ge 3$) in the projective plane, the main theorem of this paper is
\begin{theorem}\label{mainth}
When $3\le k\le 7,$ we have $cr_{N_1}(P(3k,k))= k-2$. Moreover, for $k\ge 8,$ we have $k-2\le cr_{N_1}(P(3k,k))\le k-1.$
\end{theorem}

Let us give an overview of the rest of this paper. Some basic notations are introduced in Section \ref{secpre}. Section \ref{secbou} is devoted to give the proof of Theorem \ref{mainth} by investigating the upper and lower bound of $cr_{N_1}(P(3k,k))$ independently. The lower bound of $cr_{N_1}(P(3k,k))$ is based on the result of Lemma \ref{leD}, which will be proved finally in Section \ref{secle}.

\section{Preliminaries}\label{secpre}

For $F\subseteq E$, we denote by $G\setminus F$ the graph obtained from $G$ by deleting all edges in $F$. Furthermore, we also use $F$ to denote the subgraph induced on the edge set $F$ if there is no ambiguity.

Let $k$ be an integer greater than or equal to 3, the {\it generalized Petersen graph} $P(3k,k)$  is the graph with vertex set
\begin{eqnarray*}
V(P(3k,k))=\{u_i, v_i| i=1,2,\cdots ,3k\}
\end{eqnarray*}
and edge set
\begin{eqnarray*}
E(P(3k,k))= \{u_iu_{i+1}, u_iv_i, v_iv_{k+i}| i=1,2,\cdots ,3k\}
\end{eqnarray*}
the subscripts are read modulo $3k.$

To seek the structural of $P(3k,k)$, we find it is helpful to partite the edge set $E(P(3k,k))$ into several subsets $E_i$ and $H_i$ as follows. For $1\le i\le k,$ let
\begin{eqnarray*}
E_i=\{v_iv_{k+i}, v_{k+i}v_{2k+i}, v_{2k+i}v_{i}, u_iv_i, u_{k+i}v_{k+i}, u_{2k+i}v_{2k+i}\}
\end{eqnarray*}
and let
\begin{eqnarray*}
H_i=\{u_iu_{i+1}, u_{k+i}u_{k+i+1}, u_{2k+i}u_{2k+i+1}\},
\end{eqnarray*}
the subscripts are expressed modulo $3k.$

Set $E'=\bigcup_{i=1}^{k}E_i.$ It is not difficult to see that
\begin{eqnarray}\label{eqEP}
E(P(3k,k))=E' \cup \Big(\bigcup_{i=1}^{k}H_i\Big)
\end{eqnarray}
and that $P(3k,k)\setminus E_i$ contains a subgraph homeomorphic to $P(3(k-1),(k-1))$ for each $1\leq i\leq k$.

Note that three edges $v_iv_{k+i}, v_{k+i}v_{2k+i}$ and $v_{2k+i}v_{i}$ form a 3-cycle in $E_i$, which is denoted by $EC_i$ throughout the following discussions. Formally,
\begin{eqnarray*}
EC_i=v_iv_{k+i}v_{2k+i}v_{i}.
\end{eqnarray*}

Let $D$ be a good drawing of $P(3k,k)$ in the projective plane, we define a function $f_D(H_i)$ ($1\leq i\leq k$) counting the number of crossings related to $H_i$ in $D$ as follows:
\begin{eqnarray}\label{eqfd}
f_D(H_i)=v_D(H_i,H_i)+\frac{1}{2}\sum\limits_{1\leq j\leq k, \,j\neq i}v_D(H_i,H_j).
\end{eqnarray}

In the rest discussions, we find the following definition is useful.

\begin{definition}\label{def1}
An $E'-$clean drawing of $P(3k,k)$ is a good drawing of $P(3k,k)$ in the projective plane such that $E'$ is clean.
\end{definition}

By Eqs.(\ref{eqEP}) and (\ref{eqfd}) and by counting the number of crossings in $D,$ we get
\begin{lemma}\label{lefd}
Let $D$ be an $E'-$clean drawing of $P(3k,k)$, then
$$v(D)=\sum\limits_{i=1}^{k}f_D(H_i).$$
\end{lemma}

\section{The sketch of the proof of Theorem \ref{mainth}}\label{secbou}

In the beginning of this section, the upper bound of $cr_{N_1}(P(3k,k))$ can be obtained easily.

\begin{lemma}\label{leup}
$cr_{N_1}(P(9,3))\le 1$, and $cr_{N_1}(P(3k,k))\le k-1$ for $k\ge 4.$
\end{lemma}

\noindent{\bf Proof.} Wilson's Lemma states that the crossing number of a non-planar graph in the projective plane is strictly less than its crossing number in the plane\cite{Riskin1993}. Combining this fact with the result that $cr(P(9,3))=2$ \cite{FS1986} and that $cr(P(3k,k))=k$ for $k\ge 4$ \cite{Wang2011}, the lemma follows easily.  \hfill$\Box$

\vskip 2mm
Our main efforts are made to establish the lower bound of $cr_{N_1}(P(3k,k)),$ which is based on the lemma below.

\begin{lemma}\label{leD}
For $k\ge 4,$ let $D$ be an $E'-$clean drawing of $P(3k,k)$. Then $v(D)\ge k-1.$
\end{lemma}

We postpone its proof to Section \ref{secle}. By assuming Lemma \ref{leD}, we can prove the lower bound of $cr_{N_1}(P(3k,k))$ immediately.

\begin{lemma}\label{lelow}
$cr_{N_1}(P(3k,k))\ge k-2$ for $k\ge 3.$
\end{lemma}

\begin{figure}[htbp]
\begin{minipage}[t] {0.5\linewidth}
\centering
\includegraphics[width=0.78\textwidth]{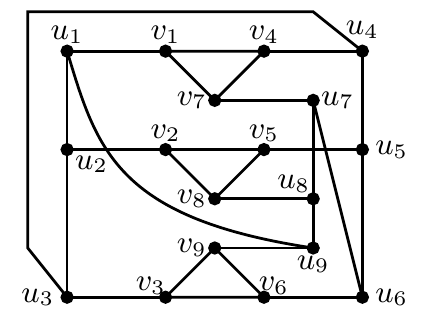}
\caption{The graph $P(9,3).$}
\label{fig1}
\end{minipage}%
\quad
\begin{minipage}[t] {0.48\linewidth}
\centering
\includegraphics[width=0.9\textwidth]{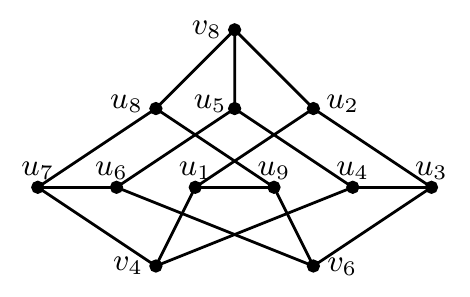}
\caption{The graph $F_{13}(12,18).$ }
\label{fig2}
\end{minipage}%
\end{figure}

\noindent{\bf Proof.} We prove the lemma by induction on $k.$ First of all, it is seen from Figure \ref{fig1} and Figure \ref{fig2} that $P(9,3)\setminus \{v_3v_9, v_2v_5, v_1v_7\}$ is a subdivision of $F_{13}(12,18)$, which is one of minimal forbidden subgraphs for the projective plane (see Appendix A in \cite{Mohar2001}), therefore, the induction basis $cr_{N_1}(P(9,3))\ge 1$ holds. Suppose that $cr_{N_1}(P(3(k-1),(k-1)))\ge k-3$ when $k\ge 4,$ consider now the graph $P(3k,k).$ Let $D$ be a good drawing of $P(3k,k)$ in $N_1.$

\vskip 2mm
{\bf Case 1.} There exists an integer $i$ $(1\le i\le k)$ such that $E_i$ is crossed in $D$.

W.l.o.g., we may assume that $i=1.$  By deleting the edges of $E_1$ in $D$, we get a good drawing $D_0$ of $P(3(k-1),(k-1))$ in $N_1$ with at least $k-3$ crossings by the induction hypothesis, thus
$$v(D)\geq v(D_0)+1\geq (k-3)+1=k-2.$$

\vskip 2mm
{\bf Case 2.} $E_i$ is clean in $D$ for every $1\le i\le k$.

Then $D$ is an $E'-$clean drawing of $P(3k,k)$. Due to Lemma \ref{leD}, $v(D)\ge k-1$ holds.

Therefore, $cr_{N_1}(P(3k,k))\ge k-2$ when $k\ge 3.$  \hfill$\Box$

\vskip 2mm

Based on the former lemmas, we can prove Theorem \ref{mainth} in the following.

\vskip 2mm

\begin{figure}[htbp]
\begin{minipage}[t] {0.49\linewidth}
\centering
\includegraphics[width=0.9\textwidth]{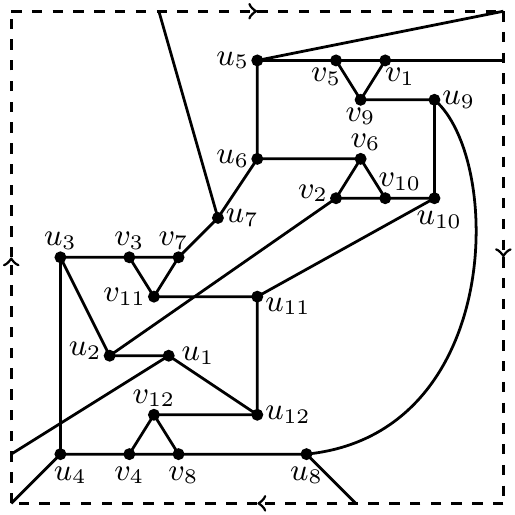}
\caption{A good drawing of $P(12,4)$ in $N_1$ with two crossings.}
\label{figp12}
\end{minipage}%
\quad
\begin{minipage}[t] {0.49\linewidth}
\centering
\includegraphics[width=0.9\textwidth]{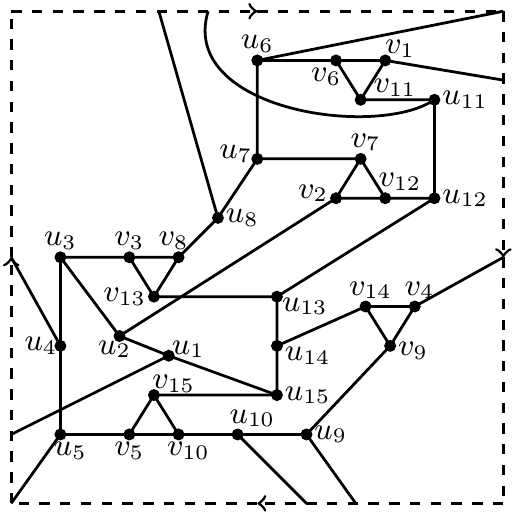}
\caption{A good drawing of $P(15,5)$ in $N_1$ with three crossings. }
\label{figp15}
\end{minipage}%
\end{figure}

\begin{figure}[htbp]
\begin{minipage}[t] {0.49\linewidth}
\centering
\includegraphics[width=0.9\textwidth]{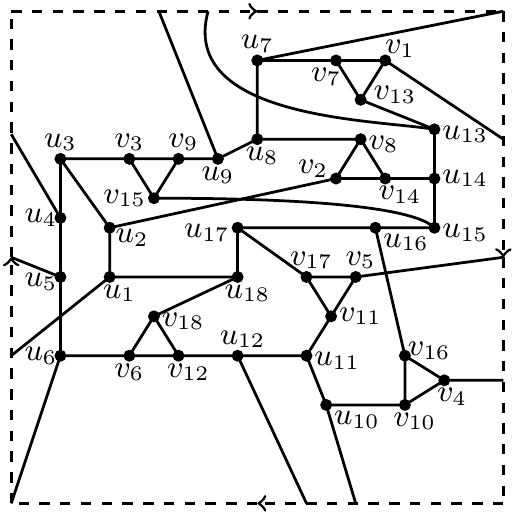}
\caption{A good drawing of $P(18,6)$  in $N_1$ with four crossings.}
\label{figp18}
\end{minipage}%
\quad
\begin{minipage}[t] {0.49\linewidth}
\centering
\includegraphics[width=0.9\textwidth]{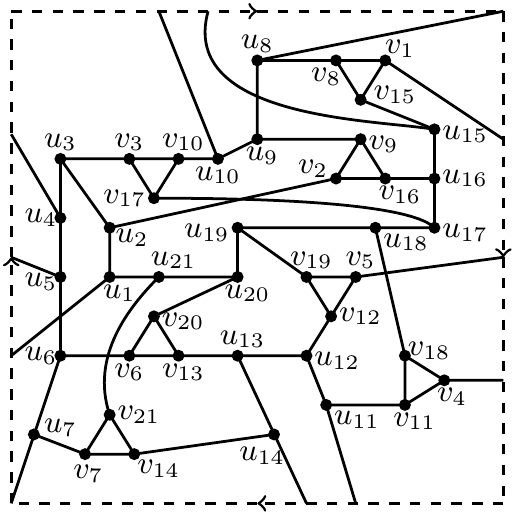}
\caption{A good drawing of $P(21,7)$ in $N_1$ with five crossings. }
\label{figp21}
\end{minipage}%
\end{figure}

\noindent{\bf The proof of Theorem \ref{mainth}.}  Figures \ref{figp12}, \ref{figp15}, \ref{figp18} and \ref{figp21} demonstrate good drawings of $P(3k,k)$ in $N_1$ with $k-2$ crossings when $4\le k\le 7,$ thus $cr_{N_1}(P(3k,k))\le k-2.$ Together with Lemma \ref{leup} and Lemma \ref{lelow}, it is confirmed that $cr_{N_1}(P(3k,k))= k-2$ for $3\le k\le 7.$

For $k\ge 8,$ we have $k-2\le cr_{N_1}(P(3k,k))\le k-1$ due to Lemma \ref{leup} and Lemma \ref{lelow}. The proof is completed. \hfill$\Box$

\section{The proof of Lemma \ref{leD}}\label{secle}

For $1\le i\le k,$ let
\begin{eqnarray}\label{eqRi}
R_i=E_i\cup H_i \cup E_{i+1}.
\end{eqnarray}

The following observation is obvious, since it is appeared several times in the rest of the paper, we state it formally.

\begin{observation}\label{obser1}
Let $D$ be an $E'-$clean drawing of $P(3k,k)$, then $E_i$ is clean for all $1\le i\le k.$
\end{observation}

\begin{figure}[htbp]
\begin{minipage}[t] {0.48\linewidth}
\centering
\includegraphics[width=0.9\textwidth]{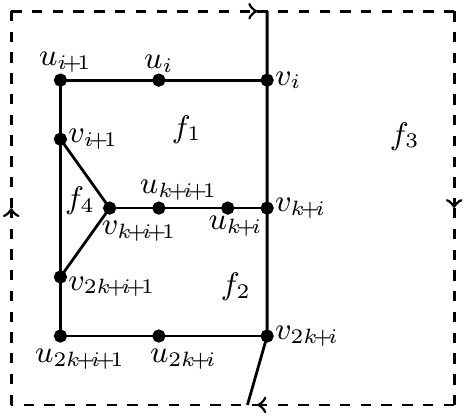}
\caption{A subdrawing of $R_i.$}
\label{fig4}
\end{minipage}%
\quad
\begin{minipage}[t] {0.48\linewidth}
\centering
\includegraphics[width=0.9\textwidth]{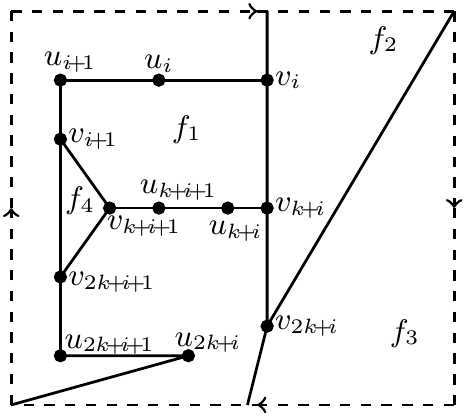}
\caption{A subdrawing of $R_i.$ }
\label{fig5}
\end{minipage}%
\end{figure}

\begin{lemma}\label{leCnon}
For $k\ge 3,$ let $D$ be an $E'-$clean drawing of $P(3k,k)$ such that one of the 3-cycles $EC_i$ and $EC_{i+1}$ is non-contractible. Then $f_D(H_i)\ge 1.$
\end{lemma}

\noindent{\bf Proof.} W.l.o.g., assume that $EC_i$ is non-contractible, for the other case, the proof is analogous. Then $EC_{i+1}$ is contractible, otherwise, two non-contractible cycles $EC_i$ and $EC_{i+1}$ will cross each other at least once in $N_1,$ a contradiction with Observation \ref{obser1}.

\begin{figure}[htbp]
\begin{minipage}[t] {0.48\linewidth}
\centering
\includegraphics[width=0.85\textwidth]{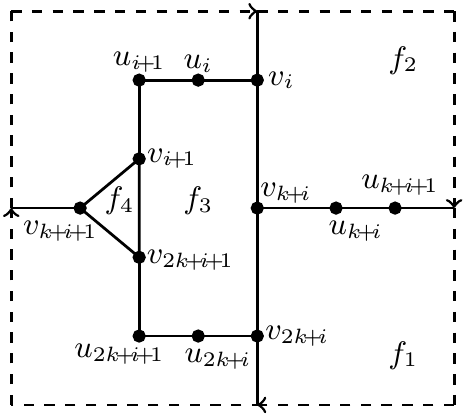}
\caption{A subdrawing of $R_i.$}
\label{fig6}
\end{minipage}%
\quad
\begin{minipage}[t] {0.48\linewidth}
\centering
\includegraphics[width=0.75\textwidth]{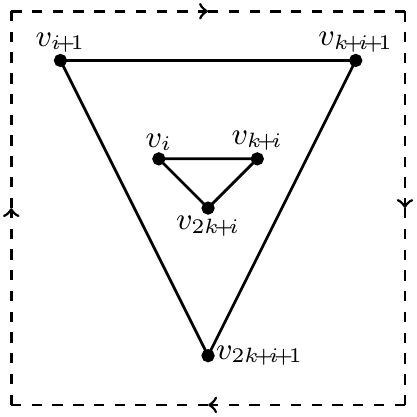}
\caption{A subdrawing of $EC_i\cup EC_{i+1}.$ }
\label{fig7}
\end{minipage}%
\end{figure}

Suppose to contrary that $f_D(H_i)< 1,$ then three edges of $H_i$ cannot cross each other by Eq.(\ref{eqfd}). By the above analyses, there are three possibilities of the subdrawing of $R_i$, see Figure \ref{fig4}, Figure \ref{fig5} and Figure \ref{fig6}. It is seen that, in each subdrawing of $R_i,$ the projective plane has been divided into four regions.

Firstly, we consider the case that the subdrawing of $R_i$ is as drawn in Figure \ref{fig4}. Consider now the subgraph $E_{i+2}$. It is asserted that all of the vertices of $E_{i+2}$ lie in the same region in the subdrawing of $R_i$ by Observation \ref{obser1}, furthermore, they cannot lie in the region labelled $f_4,$ otherwise the edge $u_{i+1}u_{i+2}$ will cross $EC_{i+1}$ at least once, a contradiction with Observation \ref{obser1}. The following three cases are investigated.

\vskip 2mm
{\bf Case 1.} The vertices of $E_{i+2}$ lie in the region labelled $f_1.$

Then the edge $u_{2k+i+1}u_{2k+i+2}$ and the path $u_{k+i+2}u_{k+i+3}\cdots u_{2k+i}$ will cross $H_i$ at least once, respectively, which implies that $f_D(H_i)\ge 1$ by Eq.(\ref{eqfd}), a contradiction.

\vskip 2mm
{\bf Case 2.} The vertices of $E_{i+2}$ lie in the region labelled $f_2.$

Then the edge $u_{i+1}u_{i+2}$ and the path $u_{2k+i+2}u_{2k+i+3}\cdots u_{3k}u_1\cdots u_{i}$ will cross $H_i$ at least once, respectively, which implies that $f_D(H_i)\ge 1$, a contradiction.

\vskip 2mm
{\bf Case 3.} The vertices of $E_{i+2}$ lie in the region labelled $f_3.$

Then $H_i$ will be crossed by the edge $u_{k+i+1}u_{k+i+2}$ and by the path $u_{i+2}u_{i+3}\cdots u_{k+i}$ at least once respectively, which yields that $f_D(H_i)\ge 1$, a contradiction.

Similar contradictions can be made if the subdrawing of $R_i$ is as drawn in Figure \ref{fig5} or Figure \ref{fig6}. \hfill$\Box$

\begin{figure}[htbp]
\begin{minipage}[t] {0.48\linewidth}
\centering
\includegraphics[width=0.75\textwidth]{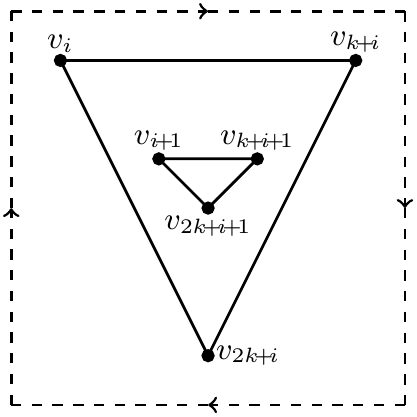}
\caption{A subdrawing of $EC_i\cup EC_{i+1}.$}
\label{fig8}
\end{minipage}%
\quad
\begin{minipage}[t] {0.48\linewidth}
\centering
\includegraphics[width=0.75\textwidth]{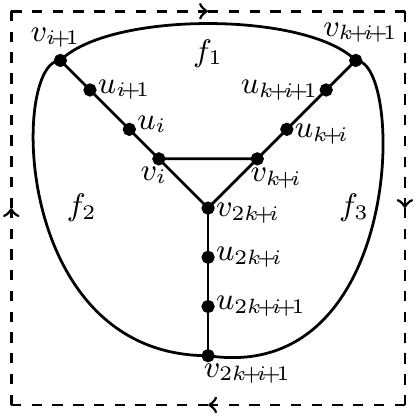}
\caption{A subdrawing of $R_i.$ }
\label{fig9}
\end{minipage}%
\end{figure}

\begin{lemma}\label{leCC}
For $k\ge 3,$ let $D$ be an $E'-$clean drawing of $P(3k,k)$. If $EC_i\cup EC_{i+1}$ is as drawn in Figure \ref{fig7} or Figure \ref{fig8}, then $f_D(H_i)\ge 1.$
\end{lemma}

\noindent{\bf Proof.} Suppose to contrary that $f_D(H_i)< 1.$  By Eq.(\ref{eqfd}), it is claimed that

\begin{claim}\label{claim1}
The edges of $H_i$ cannot cross each other in $D$.
\end{claim}

If $EC_i\cup EC_{i+1}$ is as drawn in Figure \ref{fig7}. Then $R_i$ is as shown in Figure \ref{fig9} by Observation \ref{obser1} and Claim \ref{claim1}. Consider the subgraph $E_{i+2},$ it lies in one of the regions labelled $f_1,$ $f_2$ and $f_3$  by Observation \ref{obser1} again. Assume that $E_{i+2}$ lies in the region labelled $f_1,$ then the edge $u_{2k+i+1}u_{2k+i+2}$ and the path $u_{k+i+2}u_{k+i+3}\cdots u_{2k+i}$ will cross $H_i$ at least once respectively. Hence, by Eq.(\ref{eqfd}), we have $f_D(H_i)\ge 1,$ a contradiction. Similar contradictions can be made if $E_{i+2}$ lies in the region labelled $f_2$ or $f_3.$

Using the analogous arguments, we can also show that $f_D(H_i)\ge 1$ if $EC_i\cup EC_{i+1}$ is as drawn in Figure \ref{fig8}. \hfill$\Box$

\begin{figure}[htbp]
\begin{minipage}[t] {0.48\linewidth}
\centering
\includegraphics[width=0.8\textwidth]{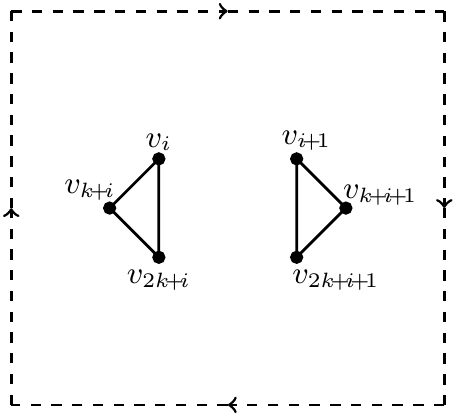}
\caption{A subdrawing of $EC_i\cup EC_{i+1}.$}
\label{fig10}
\end{minipage}%
\quad
\begin{minipage}[t] {0.48\linewidth}
\centering
\includegraphics[width=0.82\textwidth]{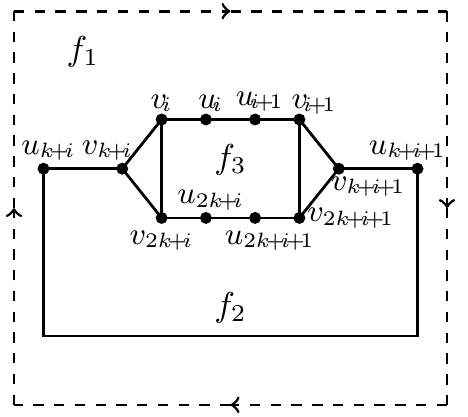}
\caption{A subdrawing of $R_i.$ }
\label{fig11}
\end{minipage}%
\end{figure}

\begin{figure}[htbp]
\begin{minipage}[t] {0.48\linewidth}
\centering
\includegraphics[width=0.85\textwidth]{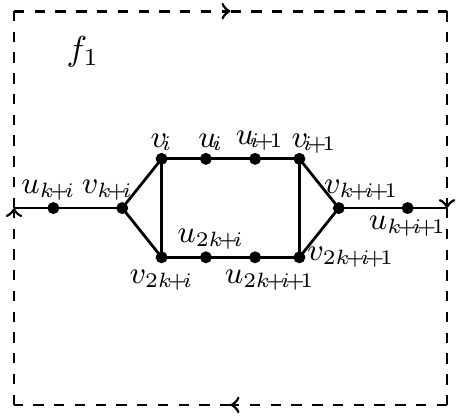}
\caption{A subdrawing of $R_i.$}
\label{fig12}
\end{minipage}%
\quad
\begin{minipage}[t] {0.48\linewidth}
\centering
\includegraphics[width=0.85\textwidth]{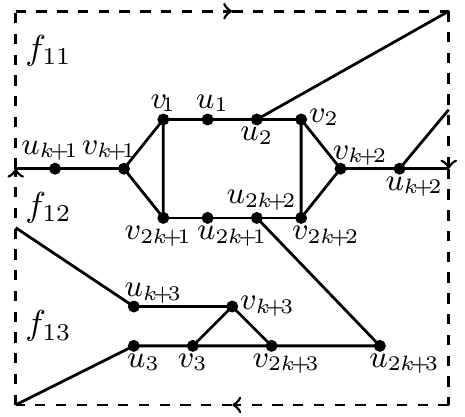}
\caption{A subdrawing of $R_1\cup H_2\cup E_3.$ }
\label{fig13}
\end{minipage}%
\end{figure}

\begin{lemma}\label{leCR1}
For $k\ge 3,$ let $D$ be an $E'-$clean drawing of $P(3k,k)$. If $EC_i\cup EC_{i+1}$ is as drawn in Figure \ref{fig10} and $f_D(H_i)< 1,$ then $R_i$ is as shown in Figure \ref{fig12}.
\end{lemma}

\noindent{\bf Proof.} First of all, we conclude that the edges of $H_i$ cannot have internal crossings in $D$ by the assumption that $f_D(H_i)< 1$. Thus, there are two possibilities of the subdrawing of $R_i$ in $D,$ see Figure \ref{fig11} and Figure \ref{fig12}.

If $R_i$ is as drawn in Figure \ref{fig11}. Observation \ref{obser1} enforces that $E_{i+2}$ lies in one of the region labelled $f_1,$ $f_2$ and $f_3$. We may assume firstly that $E_{i+2}$ lies in the region labelled $f_1,$ for other cases the proof is the same. Under this circumstance, the edge $u_{2k+i+1}u_{2k+i+2}$ and the path $u_{k+i+2}u_{k+i+3}\cdots u_{2k+i}$ will cross $H_i$ at least once respectively, therefore, $f_D(H_i)\ge 1$ by Eq.(\ref{eqfd}), a contradiction.

Thus, $R_i$ is as shown in Figure \ref{fig12}. \hfill$\Box$

\vskip 2mm
Combining Lemma \ref{leCC} with Lemma \ref{leCR1}, we have

\begin{lemma}\label{leCR2}
For $k\ge 3,$ let $D$ be an $E'-$clean drawing of $P(3k,k)$ such that both $EC_i$ and $EC_{i+1}$ are contractible cycles. If $R_i$ is not as drawn in Figure \ref{fig12}, then $f_D(H_i)\ge 1.$
\end{lemma}

\noindent{\bf Proof.} Since both $EC_i$ and $EC_{i+1}$ are contractible cycles, there are three possibilities of the subdrawing of $EC_i\cup EC_{i+1}$, see Figure \ref{fig7}, Figure \ref{fig8} and Figure \ref{fig10}. Lemma \ref{leCC} implies that $f_D(H_i)\ge 1$ if $EC_i\cup EC_{i+1}$ is as drawn in Figure \ref{fig7} or Figure \ref{fig8}. Furthermore, Lemma \ref{leCR1} implies that $f_D(H_i)\ge 1$ if $R_i$ is not as drawn in Figure \ref{fig12}. \hfill $\Box$

\vskip 2mm
Now, we are ready to prove Lemma \ref{leD}.
\vskip 2mm

\noindent{\bf The proof of Lemma \ref{leD}.} Suppose to contrary that
\begin{eqnarray}\label{eqvD}
v(D)<k-1.
\end{eqnarray}
By Lemma \ref{lefd}, there exists an integer $i$ such that $f_D(H_i)<1,$ otherwise,
\begin{eqnarray*}
v(D)=\sum\limits_{i=1}^{k}f_D(H_i)\ge k,
\end{eqnarray*}
a contradiction with Eq.(\ref{eqvD}). The following two cases are discussed: Case 1. there exists an integer $i$ such that $f_D(H_i)=0$ and Case 2. $f_D(H_i)>0$ for every $1\le i\le k$ and there exists  an integer $i$ such that $f_D(H_i)=\frac{1}{2}.$

\vskip 2mm
{\bf Case 1.} There exists an integer $i$ such that $f_D(H_i)=0.$

W.l.o.g., let $f_D(H_1)=0.$ Thus, we can get that

\begin{claim}\label{claim2}
$R_1$ is clean in $D.$
\end{claim}

Furthermore, there exists another integer $j$ $(j\ne 1)$ satisfying $f_D(H_j)<1,$ otherwise,
\begin{eqnarray*}
v(D)=\sum\limits_{i=2}^{k}f_D(H_i)\ge k-1,
\end{eqnarray*}
a contradiction with Eq.(\ref{eqvD}).

\vskip 2mm
{\bf Subcase 1.1.} $j=2$ or $j=k$.

By symmetry, we only need to consider the case that $j=2$. It follows from Lemma \ref{leCnon} that both $EC_1$ and $EC_2$ are contractible cycles since $f_D(H_1)=0.$ Moreover, Lemma \ref{leCR2} enforces that the subdrawing of $R_1$ is as shown in Figure \ref{fig12} by replacing all the indices $i$ by 1.

Now, we consider the subgraph $E_{3}$. All of the vertices of $E_{3}$ lie in the region labeled $f_1$, on whose boundary lies the vertices $u_2, u_{k+2}$ and $u_{2k+2},$ otherwise, we have $v_D(H_2,R_1)\ge 1$, contradicting Claim \ref{claim2}. Moreover, the edges of $H_2$ cannot have internal crossings by Eq.(\ref{eqfd}). All these arguments confirm that the only possibility of the subdrawing of $R_1\cup H_2\cup E_3$ is as shown in Figure \ref{fig13}. Note that the region labeled by $f_1$ in Figure \ref{fig12} has been separated into three regions, which are labelled by $f_{11}$, $f_{12}$ and $f_{13}$.

Next, we consider the subgraph $E_k$ (note that $k\ge 4$ is crucial here for $E_k$ being not equal to $E_3$). By Observation \ref{obser1}, $E_k$ lies in one of $f_{11}$, $f_{12}$ and $f_{13}$. We may assume that $E_k$ lies in $f_{11}$, for other cases the proof is the same. Under this circumstance, $H_2$ will be crossed by the edge $u_{2k}u_{2k+1}$ and by the path $u_{k+3}u_{k+4}\cdots u_{2k}$ at least once respectively, therefore, $f_D(H_2)\ge 1$ by Eq.(\ref{eqfd}), a contradiction.

\vskip 2mm
{\bf Subcase 1.2.} $j\notin\{2,k\}.$

Since $f_D(H_1)=0$ and $f_D(H_j)<1,$ all of the cycles $EC_1,$ $EC_2,$ $EC_j$ and $EC_{j+1}$ are contractible by Lemma \ref{leCnon}. By Lemma \ref{leCR2}, $R_1$ (resp. $R_j$) is as drawn in Figure \ref{fig12} by replacing all the indices $i$ by 1 (resp. $j$).

Notice that both $u_{k+1}v_{k+1}v_1u_1u_2v_2v_{k+2}u_{k+2}u_{k+1}$ and $u_{k+j}v_{k+j}v_ju_ju_{j+1}v_{j+1}v_{k+j+1}$ $u_{k+j+1}u_{k+j}$ are non-contractible curves, therefore, they must cross each other in $N_1$, which yields that $v_D(H_1,H_j)\ge 1$ since $E'$ is clean in $D$, a contradiction with $f_D(H_1)=0$.

\vskip 2mm
{\bf Case 2.} $f_D(H_i)>0$ for every $1\le i\le k$ and there exists an integer $i$ such that $f_D(H_i)=\frac{1}{2}.$

W.l.o.g., let $f_D(H_1)=\frac{1}{2}.$ There exists an integer $j$ ($j\notin \{2,k\}$) such that $f_D(H_j)=\frac{1}{2},$ otherwise
\begin{eqnarray*}
v(D)=\sum\limits_{i=1}^{k}f_D(H_i)\ge 3\times\frac{1}{2}+\sum\limits_{i=4}^{k}f_D(H_i)\ge \frac{3}{2}+(k-3),
\end{eqnarray*}
and thus $v(D)\ge k-1$ since $v(D)$ is an integer, a contradiction with Eq.(\ref{eqvD}).

\begin{figure}[htbp]
\centering
\includegraphics[width=0.42\textwidth]{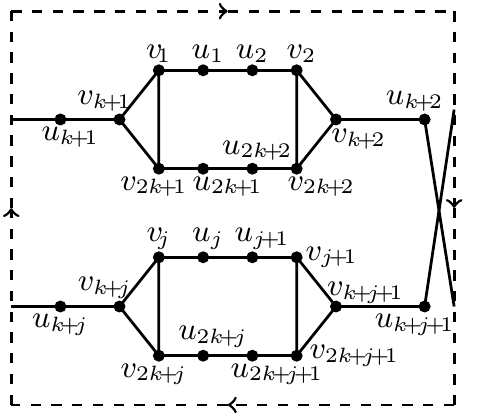}
\caption{A subdrawing of $R_1\cup R_j.$ }
\label{fig15}
\end{figure}

Lemma \ref{leCnon} tells that, all of the cycles $EC_1,$ $EC_2,$ $EC_j$ and $EC_{j+1}$ are contractible. By Lemma \ref{leCR2}, $R_1$ (resp. $R_j$) is as drawn in Figure \ref{fig12} by replacing all the indices $i$ by 1 (resp. $j$). Note that both $u_{k+1}v_{k+1}v_1u_1u_2v_2v_{k+2}u_{k+2}u_{k+1}$ and $u_{k+j}v_{k+j}v_ju_ju_{j+1}v_{j+1}v_{k+j+1}$$u_{k+j+1}u_{k+j}$ are non-contractible curves, then
$$v_D(H_1,H_j)=1$$
since $f_D(H_1)=f_D(H_j)=\frac{1}{2}$, furthermore, the crossed edge of $H_1$ (resp. $H_j$) is $u_{k+1}u_{k+2}$ (resp. $u_{k+j}u_{k+j+1}$), see Figure \ref{fig15}. It is seen that two vertices $u_2$ and $u_j$ don't lie on the boundary of a same region, thus the path $u_2u_3\cdots u_j$ will cross $H_1$ or $H_j$ by Observation \ref{obser1}, a contradiction with the assumption that $f_D(H_1)=f_D(H_j)=\frac{1}{2}.$

All of the above contradictions enforce that $v(D)\ge k-1$ if $D$ is an $E'-$clean drawing of $P(3k,k)$.  \hfill$\Box$

\section{Conclusions}

This paper studies the crossing number of the generalized Petersen graph $P(3k,k)$ in the projective plane. However, when $k\ge 8$, the exact value of $cr_{N_1}(P(3k,k))$ remains open. Does $cr_{N_1}(P(3k,k))=k-2$ or $cr_{N_1}(P(3k,k))=k-1$? From the former proof, we know the problem of "Deciding the exact value of $cr_{N_1}(P(3k,k))$" is highly related to the problem of "Whether there exists a good drawing of $P(3k,k)$ in the projective plane with $k-2$ crossings or not?" If the answer to the latter problem is "Yes", then we conclude that $cr_{N_1}(P(3k,k))=k-2$, otherwise, we have $cr_{N_1}(P(3k,k))=k-1$.


\end{document}